\documentclass[11pt]{amsart}
\usepackage{amssymb}
\usepackage{amsmath}
\usepackage[active]{srcltx}
\usepackage{t1enc}
\usepackage[latin2]{inputenc}
\usepackage{verbatim}
\usepackage{amsmath,amsfonts,amssymb,amsthm}
\usepackage[mathcal]{eucal}
\usepackage{enumerate}
\usepackage[centertags]{amsmath}
\usepackage{graphics}

\newtheorem{theorem}{Theorem}

\newtheorem{lemma}{Lemma}
\newtheorem{remark}{Remark}

\newtheorem{corollary}{Corollary}

\begin{document}
\author{David Baramidze,  Lars-Erik Persson, Harpal Singh and George Tephnadze}
\title[partial sums ]{Some new weak $(H_p-L_p)$ type inequality for weighted maximal operators of Walsh-Fourier series}

\address{D. Baramidze, The University of Georgia, School of science and technology, 77a Merab Kostava St, Tbilisi 0128, Georgia and Department of Computer Science and Computational Engineering, UiT - The Arctic University of Norway, P.O. Box 385, N-8505, Narvik, Norway.}
\email{ datobaramidze20@gmail.com \ \ \ davit.baramidze@ug.edu.ge   }
\address{L.-E. Persson, UiT The Arctic University of Norway, P.O. Box 385, N-8505, Narvik, Norway and Department of Mathematics and Computer Science, Karlstad University, 65188 Karlstad, Sweden.}
\email{larserik6pers@gmail.com}
\address{H. Singh, Department of Computer Science and Computational Engineering, UiT - The Arctic University of Norway, P.O. Box 385, N-8505, Narvik, Norway.}
\email{harpal.singh@uit.no }
\address {G. Tephnadze, The University of Georgia, School of Science and Technology, 77a Merab Kostava St, Tbilisi, 0128, Georgia.}
\email{g.tephnadze@ug.edu.ge}

\thanks{The research was supported by Shota Rustaveli National Science Foundation grant no. PHDF-21-1702.}
\maketitle

\begin{abstract}
In this paper we  introduce some new weighted  maximal operators of the partial sums of the Walsh-Fourier series. We prove that for some "optimal" weights these new operators indeed are bounded  from the martingale Hardy space $H_{p}$ to the Lebesgue space $\text{weak}-L_{p},$  for $0<p<1.$ Moreover, we also prove sharpness of this result. 

\end{abstract}

\date{}



\section{INTRODUCTION}

\bigskip  All symbols used in this introduction can be found in Section 2.

It is well-known that  the Walsh system does
not form a basis in the space $L_1(G)$ (see e.g. \cite{1}). Moreover, there exists a function in the dyadic Hardy space $H_{1}(G),$ such that the partial sums of $f$ are not bounded in the  $L_{1}$-norm.
Uniform and pointwise convergence and some approximation properties of
partial sums in $L_{1}(G)$ norms were investigated by  Avdispahić and Memić \cite{am},  Gát, Goginava and Tkebuchava \cite{ggt,gt}, Nagy \cite{na}, Onneweer \cite{10} and Persson, Schipp,  Tephnadze and Weisz \cite{PSTW}. Fine \cite{fi} obtained sufficient conditions for the uniform convergence
which are completely analogous to the Dini-Lipschits conditions. Guličev \cite{9} estimated the rate of uniform convergence of a Walsh-Fourier series by using Lebesgue constants and modulus of continuity. These problems for
Vilenkin groups were investigated by Blatota, Nagy, Persson and Tephnadze \cite{BNPT} (see also \cite{2,BNT9,BNT10}), Fridli \cite{4} and Gát  \cite{5}. 

To study convergence of subsequences of  Fejér means and their restricted  maximal operators on the martingale Hardy spaces $H_p(G)$ for $0<p\leq 1/2,$ the central role is played by the fact that any natural number $n\in \mathbb{N}$ can be uniquely expressed as
\begin{equation}\label{star01}
n=\sum_{k=0}^{\infty }n_{j}2^{j},  \ \ n_{j}\in Z_{2} ~(j\in \mathbb{N}), 
\end{equation} 
where only a finite numbers of $n_{j}$ differ from zero
and their important characters  $\left[ n\right],$ $\left\vert n\right\vert,$ $\rho\left( n\right)$  and $V(n)$ are defined by
\begin{equation}\label{star02}
\left[ n\right] :=\min \{j\in \mathbb{N},n_{j}\neq 0\}, \  
\left\vert n\right\vert :=\max \{j\in \mathbb{N},n_{j}\neq 0\}, \  
\rho\left( n\right) =\left\vert n\right\vert -\left[ n\right] 
\end{equation}
and
\begin{equation*}
V\left( n\right): =n_{0}+\overset{\infty }{\underset{k=1}{\sum }}\left|
n_{k}-n_{k-1}\right|, \text{ \ for \
	all \ \ }n\in \mathbb{N}
\end{equation*}

In particular, (see \cite{BPT1}, \cite{luk} and  \cite{sws}) 
\begin{equation*}
\frac{V\left( n\right) }{8}\leq \Vert D_n\Vert_1\leq V\left( n\right)
\end{equation*}
from which it follows that, for any $F\in L_1(G),$ there exists an absolute constant $c$ such that the following inequality holds:

\begin{equation*}
\left\Vert S_n F\right\Vert_1\leq c{V\left( n\right)
}\left\Vert F\right\Vert_1.
\end{equation*}
Moreover, for any $f\in H_1(G)$  (see \cite{tep2})
\begin{equation*}
\left\Vert S_{n}F\right\Vert _{H_{1}}\leq c{V\left( n\right)}\left\Vert F\right\Vert _{H_{1}}.
\end{equation*}

For $0<p<1$ in  \cite{tep0,tep1}  the weighted  maximal operator $\overset{\sim }{S }^{*,p},$ defined by
\begin{equation}\label{max0}
\overset{\sim }{S }^{*,p}F:=\sup_{n\in\mathbb{N}}\frac{\left|S _{n}F\right|} {\left( n+1\right)^{1/p-1} }
\end{equation}
was investigated and it was proved that the following inequalities hold:
\begin{equation*}
\left\Vert \overset{\sim }{S }^{*}F\right\Vert_{p}\leq c_{p}\left\Vert F\right\Vert _{H_{p}}
\end{equation*}
Moreover, it was also proved that the rate of the sequence $\{\left(n+1\right)^{ 1/p-1}\}$ given in the denominator of \eqref{max0} can not be improved.

In \cite{tep2} and \cite{tep3} (see also \cite{Bara3}) it was proved that if $F\in H_{p}(G),$ then there exists an absolute constant  $c_{p},$ depending only on $p,$ such that
\begin{equation*}
	\text{ }\left\Vert S_{n}F\right\Vert _{H_{p}}\leq c_{p}2^{\rho\left( n\right)\left( 1/p-1\right) }\left\Vert F\right\Vert _{H_{p}},
\end{equation*}	
which implies that

\begin{equation*}
\left\Vert \frac{S_{n}F}{2^{\rho \left(n\right) \left( 1/p-1\right) }}\right\Vert_{p}\leq c_{p}\left\Vert F\right\Vert _{H_{p}}.
\end{equation*}
Moreover, if $0<p<1,$ $\left\{ n_{k}:\text{ }k\geq 0\right\} $ is any increasing sequence of positive integers such that 
\begin{equation*}
\sup_{k\in \mathbb{N}}\rho\left( n_{k}\right) =\infty  
\end{equation*}
and  
$\Phi :\mathbb{N}_{+}\rightarrow \lbrack 1,\infty )$ 
is any nondecreasing function, satisfying the condition 
\begin{equation*}	
\overline{\underset{k\rightarrow \infty }{\lim }}\frac{2^{\rho\left(
n_{k}\right) \left( 1/p-1\right) }}{\Phi \left( n_{k}\right) }=\infty,
\end{equation*}
then there exists a martingale $F\in H_{p}(G),$ such that
\begin{equation*}
\underset{k\in \mathbb{N}}{\sup }\left\Vert \frac{S_{n_{k}}F}{\Phi \left(n_{k}\right) }\right\Vert_{\text{weak}-L_p}=\infty .
	\end{equation*}

Convergence and summability of partial sums and general summability methods of Walsh-Fourier series and boundedness of their maximal operators can be found in Baramidze \cite{bar1,bar2,bar3,bar5,BGN1}, Blahota  \cite{BT1,BBTT}, Nagy  \cite{NT1,NT2,NT3,NT4}, Persson, Tephnadze and Weisz \cite{NPTW,PSTW},  Persson, Tephnadze and Wall \cite{PTTW1,PTW2,PTW3,PTW4}, Tephnadze \cite{AT,MST,tep4,tep5,tep6,tep7,tep8,tep9,tep10}.

In this paper we prove that the weighted maximal operator  of the partial sums of the Walsh-Fourier defined by 
\begin{equation*} 
\sup_{n\in \mathbb{N}}\frac{\left\vert S_{n}F\right\vert}{2^{\rho \left(n\right) \left( 1/p-1\right)}}
\end{equation*}
is bounded from the martingale Hardy space $H_p(G)$ to the space $weak-L_p(G),$  for $0<p<1.$ We also prove the sharpness of this result (see Theorem 2). As a consequence we obtain both some new and well-known results.

This paper is organized as follows: In order not to disturb our discussions later on some preliminaries are presented in Section 2. The main results and some of its consequences can be found in Section 3.  The detailed proofs of the main results are given in Section 4. Finally, Section 5 contains some final remarks related to the new book \cite{tepbook} and two related conjectures.

\section{Preliminaries}

\bigskip Let $\mathbb{N}_{+}$ denote the set of the positive integers, $%
\mathbb{N}:=\mathbb{N}_{+}\cup \{0\}.$ Denote by $Z_{2}$ the discrete cyclic
group of order $2$, that is $Z_{2}:=\{0,1\},$ where the group operation is the modulo $2$ addition and every subset is open. The Haar measure on $Z_{2}$ is given so that the measure of a singleton is $1/2$.

Define the group $G$ as the complete direct product of the group $Z_{2},$
with the product of the discrete topologies of $Z_{2}$`s. The elements of $G$
are represented by sequences 
$$x:=(x_{0},x_{1},...,x_{j},...), \ \ \ \text{ where } \ \ \
x_{k}=0\vee 1.$$

It is easy to give a base for the neighborhood of $x\in G:$ 
\begin{equation*}
I_{0}\left( x\right) :=G,\text{ \ }I_{n}(x):=\{y\in
G:y_{0}=x_{0},...,y_{n-1}=x_{n-1}\}\text{ }(n\in \mathbb{N}).
\end{equation*}

Denote 
$I_{n}:=I_{n}\left( 0\right) , \ \  \overline{I_{n}}:=G\backslash 
I_{n}$
and
$$  e_{n}:=\left( 0,...,0,x_{n}=1,0,...\right) \in G, \ \  \text{ for } \ \  n\in \mathbb{N}.$$ 
Then it is easy to prove that
\begin{equation}\label{2}
\overline{I_{M}}=\overset{M-1}{\underset{s=0}{\bigcup }}I_{s}\backslash
I_{s+1}.  
\end{equation}

The norms (or quasi-norms) of the Lebesgue space $L_{p}(G)$ and the weak Lebesgue space $L_{p,\infty }\left(
G\right) ,$ $\left( 0<p<\infty \right) $ are, respectively, defined by 
\begin{equation*}
\left\Vert f\right\Vert _{p}^{p}:=\int_{G}\left\vert f\right\vert ^{p}d\mu
\ \ \ \text{
and} \ \ \
\left\Vert f\right\Vert _{weak-L_{p}}^{p}:=\sup_{\lambda
	>0}\lambda ^{p}\mu \left( f>\lambda \right) .
\end{equation*}

The $k$-th Rademacher function $r_{k}\left( x\right)$ is defined by
\begin{equation*}
r_{k}\left( x\right) :=\left( -1\right) ^{x_{k}}\text{\qquad }\left( \text{ }x\in G,\text{ }k\in \mathbb{N}\right) .
\end{equation*}

Now, define the Walsh system $w:=(w_{n}:n\in \mathbb{N})$ on $G$ by
\begin{equation*}
w_{n}(x):=\overset{\infty }{\underset{k=0}{\Pi }}r_{k}^{n_{k}}\left(
x\right) =r_{\left\vert n\right\vert }\left( x\right) \left( -1\right) ^{%
\underset{k=0}{\overset{\left\vert n\right\vert -1}{\sum }}n_{k}x_{k}}\text{%
\qquad }\left( n\in \mathbb{N}\right) .
\end{equation*}

The Walsh system is orthonormal and complete in $L_{2}\left( G\right) $ (see
e.g. \cite{sws}).

If $f\in L_{1}\left( G\right) $ we can establish the Fourier coefficients,
the partial sums of the Fourier series, the Dirichlet kernels with respect
to the Walsh system in the usual manner:$\qquad $
\begin{eqnarray*}
\widehat{f}\left(k\right)&:=&\int_{G}fw_{k}d\mu \,\,\,\,\left( k\in \mathbb{N}\right), \\ \notag
\end{eqnarray*}
\begin{eqnarray*}
S_{n}f&:=&\sum_{k=0}^{n-1}\widehat{f}\left( k\right) w_{k}, \,\,\,\left( n\in \mathbb{N}_{+}\right)\\ \notag
D_{n}&:=&\sum_{k=0}^{n-1}w_{k\text{ }}\,\,\,\left( n\in \mathbb{N}_{+}\right). \notag
\end{eqnarray*}

Recall that (see \cite{tepbook} and \cite{sws})
\begin{equation}\label{1dn}
D_{2^{n}}\left( x\right) =\left\{ 
\begin{array}{ll}
2^{n}, & \,\text{if\thinspace \thinspace \thinspace }x\in I_{n} \\ 
0, & \text{if}\,\, \ x\notin I_{n}
\end{array}
\right.  
\end{equation}
and
\begin{equation}\label{2dn}
D_{n}=w_{n}\overset{\infty }{\underset{k=0}{\sum }}n_{k}r_{k}D_{2^{k}}=w_{n}%
\overset{\infty }{\underset{k=0}{\sum }}n_{k}\left(
D_{2^{k+1}}-D_{2^{k}}\right) ,\text{ for \ }n=\overset{\infty }{\underset{i=0%
}{\sum }}n_{i}2^{i}.  
\end{equation}
Moreover, we have the following lower estimate (see \cite{tepbook}):
\begin{lemma}\label{lemma1}
Let $n\in \mathbb{N}$ and $\left[ n\right] \neq\left\vert n\right\vert .$  Then
\begin{equation*}
\left\vert D_n(x)\right\vert =\left\vert D_{n-2^{\left\vert n\right\vert}}(x)\right\vert \geq \frac{2^{\left[ n\right] }}{4}, \ \ \ \text{for } \ \ \ x\in I_{\left[ n\right]+1}\left(e_{\left[ n\right]}\right).
\end{equation*}  
\end{lemma}

The $\sigma $-algebra generated by the intervals $\left\{ I_{n}\left(
x\right) :x\in G\right\} $ will be denoted by $\zeta _{n}\left( n\in \mathbb{%
N}\right).$ 

Denote by $F=\left( F_{n},n\in \mathbb{N}\right) $ the
martingale with respect to $\digamma _{n}$ $\left( n\in \mathbb{N}\right) $
(see e.g. \cite{We1}).

The maximal function $F^{\ast }$ of a martingale $F$ is defined by
\begin{equation*}
F^{\ast }:=\sup_{n\in \mathbb{N}}\left\vert F_{n}\right\vert .
\end{equation*}

In the case $f\in L_{1}\left( G\right),$ the maximal function $f^{\ast }$ is
given by
\begin{equation*}
f^{\ast }\left( x\right) :=\sup\limits_{n\in \mathbb{N}}\frac{1}{\mu \left(
	I_{n}\left( x\right) \right) }\left\vert \int\limits_{I_{n}\left( x\right)
}f\left( u\right) d\mu \left( u\right) \right\vert .
\end{equation*}

For $0<p<\infty $ the Hardy martingale spaces $H_{p}\left( G\right) $
consists of all martingales for which
\begin{equation*}
\left\Vert F\right\Vert _{H_{p}}:=\left\Vert F^{\ast }\right\Vert
_{p}<\infty .
\end{equation*}

It is easy to check that for every martingale $F=\left( F_{n},n\in \mathbb{N}
\right) $ and every $k\in \mathbb{N}$ the limit
\begin{equation*}
\widehat{F}\left( k\right) :=\lim_{n\rightarrow \infty }\int_{G}F_{n}\left(
x\right) w_{k}\left( x\right) d\mu \left( x\right)  
\end{equation*}
exists and it is called the $k$-th Walsh-Fourier coefficients of $F.$

If $F:=$ $\left( S_{2^n}f:n\in \mathbb{N}\right) $ is a regular martingale,
generated by $f\in L_{1}\left( G\right) ,$ then  (see e.g. \cite{tepbook}, \cite{S} and \cite{We1})
$$\widehat{F}\left(k\right) =\widehat{f}\left( k\right) , \  k\in \mathbb{N}.$$

A bounded measurable function $a$ is called $p$-atom, if there exists a dyadic
interval $I,$ such that 
\begin{equation*}
\int_{I}ad\mu =0,\text{ \ \ }\left\Vert a\right\Vert _{\infty }\leq \mu
\left( I\right) ^{-1/p},\text{ \ \ supp}\left( a\right) \subset I.
\end{equation*}

The dyadic Hardy martingale spaces $H_{p}(G)$ for $0<p\leq 1$ have an atomic
characterization. Namely, the following holds (see \cite{tepbook}, \cite{We1,We3}):

\begin{lemma}\label{W1}
A martingale $F=\left( F_{n},n\in \mathbb{N}\right) $ belongs to $H_{p}(G) \ \left( 0<p\leq 1\right) $ if and only if there exists a sequence $%
\left( a_{k},\text{ }k\in \mathbb{N}\right) $ of p-atoms and a sequence $%
\left( \mu _{k},k\in \mathbb{N}\right) $ of  real numbers such that for
every $n\in \mathbb{N},$

\begin{equation}
\qquad \sum_{k=0}^{\infty }\mu _{k}S_{2^{n}}a_{k}=F_{n},\text{ \ \ \ }%
\sum_{k=0}^{\infty }\left\vert \mu _{k}\right\vert ^{p}<\infty ,  \label{2A}
\end{equation}
Moreover, 
$
\left\Vert F\right\Vert _{H_{p}}\backsim \inf \left( \sum_{k=0}^{\infty
}\left\vert \mu _{k}\right\vert ^{p}\right) ^{1/p},
$
where the infimum is taken over all decomposition of $F$ of the form (\ref%
{2A}).
\end{lemma}

\section{The Main Results with Applications}

Our first main result reads:
\begin{theorem}\label{theorem1}
Let $0<p<1,$ $f\in {{H}_{p}}\left(G \right)$, $n$ be defined by \eqref{star01} and $\rho\left( n\right)$ be defined by \eqref{star02}. Then the weighted maximal operator $\widetilde{S }^{\ast ,\nabla },$ defined by
	\begin{equation}\label{maxoperator}
	\widetilde{S }^{\ast ,\nabla }F=\underset{n\in \mathbb{N}}{\sup }
	\frac{\left\vert S_{n}F\right\vert}{2^{\rho\left( n\right)\left( 1/p-1\right)}},
	\end{equation}
	is bounded from the martingale Hardy space ${{H}_{p}}(G)$ to the space $\text{weak}-L_p(G).$
\end{theorem}

Our second main result shows that Theorem \ref{theorem1} can not be improved in general because it is sharp in some special senses:
\begin{theorem}\label{theorem2}
a) Let $0<p<1,$ $n$ be defined by \eqref{star01},  $\rho\left( n\right)$ be defined by \eqref{star02} and $\widetilde{S }^{\ast ,\nabla }$ is defined by \eqref{maxoperator}. Then there exists a sequence $\{f_n, n\in \mathbb{N}\}$  of $p$-atoms, such that
\begin{equation*}
\sup_{n\in \mathbb{N}}\frac{\left\Vert \widetilde{S }^{\ast ,\nabla}f_n\right\Vert_{p}}{\left\Vert f_n \right\Vert_{H_p(G)}}=\infty .
\end{equation*}

b) Let $0<p<1$, $n$ be defined by \eqref{star01} and $\rho\left( n\right)$ be defined by \eqref{star02}. If $\varphi :\mathbb{N}\rightarrow \lbrack 1,$ $\infty )$ is a nondecreasing function, satisfying the condition
	\begin{equation} \label{6aa}
	\overline{\lim_{n\rightarrow \infty }}\frac{2^{\rho(n){(1/p-1)}}}{\varphi \left( n\right) }=\infty,
	\end{equation}
	then there exists a sequence $\{f_n, n\in \mathbb{N}\}$  of $p$-atoms, such that
	\begin{equation*}
	\sup_{n\in \mathbb{N}}\frac{\left\Vert \sup_{k\in \mathbb{N}}\frac{\vert S_kf_n\vert}{\varphi \left( k\right)}\right\Vert_{\text{weak}-L_p(G)}}{\left\Vert f_n \right\Vert_{H_p(G)}}=\infty .
	\end{equation*}
\end{theorem}

Theorem \ref{theorem1} implies the following result of Weisz \cite{We3} (see also \cite{We1}):
\begin{corollary}
	Let $0<p<1$ and $f\in {{H}_{p}}\left(G \right)$.
	Then the  maximal operator $S ^{\ast ,\triangle }$ defined by	
	\begin{equation*}
S ^{\ast ,\triangle }F:=	\underset{n\in \mathbb{N}}{\sup }{\left\vert S_{2^{n}}F\right\vert}
	\end{equation*}
	is bounded from the Hardy space ${{H}_{p}(G)}$ to the Lebesgue space $\text{weak}-L_p(G)$ (and, thus, to the Lebesgue space $L_p(G)$).
\end{corollary}

Moreover, Theorems 1 and 2 imply the following results  (see  \cite{tepbook}):
\begin{corollary}\label{cor2}
	Let $0<p<1$ and $f\in {{H}_{p}}\left(G \right)$.
	Then the  maximal operator $\widetilde{S }^{\ast ,\nabla },$ defined by	
	\begin{equation*}
	S^{\ast,\nabla}F:=\underset{k\in \mathbb{N}}{\sup}
	\left\vert S_{n_{k}}F\right\vert
	\end{equation*}
	is bounded from the Hardy space ${{H}_{p}(G)}$ to the Lebesgue space $\text{weak}-L_p(G)$ if and only if condition 
\begin{equation*}
\sup_{k\in \mathbb{N}}\rho\left( n_{k}\right)<c<\infty  
\end{equation*}	
	 is fulfilled.
\end{corollary}

\begin{remark}
The statement in Corollary \ref{cor2} holds also if the space $weak-L_p(G)$ is replaced by $L_p(G)$.
\end{remark}

\begin{corollary}
	a) 	Let $0<p<1$ and $f\in {{H}_{p}}\left(G \right)$.
	Then the weighted maximal operator defined by	
	\begin{equation*}
	\underset{n\in \mathbb{N}}{\sup }
	\frac{\left\vert S_{2^{n}+2^{n/2}}F\right\vert}{2^{\frac{n}{2}\left( 1/p-1\right)}}
	\end{equation*}
	is bounded from the Hardy space ${{H}_{p}}(G)$ to the Lebesgue space $\text{weak}-L_p(G).$
	
	b) (Sharpness) Let $\varphi :\mathbb{N}\rightarrow \lbrack 1,$ $\infty )$ be a nondecreasing function, satisfying the condition	
	\begin{equation*}
	\overline{\lim_{n\rightarrow \infty }}\frac{2^{\frac{n}{2}\left( 1/p-1\right)}}{\varphi \left( n\right)}=\infty. 
	\end{equation*}
	Then, there exists sequence $\{f_n, n\in \mathbb{N}\}$  of $p$-atoms, such that
	\begin{equation*}
	\sup_{n\in \mathbb{N}}\frac{\left\Vert \frac{S_{2^{n}+2^{n/2}}f_n}{\varphi \left( 2^{n}+2^{n/2}\right)} \right\Vert_{\text{weak}-L_p(G)}}{\left\Vert f_n \right\Vert_{H_p(G)}}=\infty .
	\end{equation*}
\end{corollary}

\begin{corollary}
	a) 	Let $0<p<1$ and $f\in {{H}_{p}}\left(G \right).$
	Then the weighted maximal operator defined by
	\begin{equation*}
	\underset{n\in \mathbb{N}}{\sup }
	\frac{\left\vert S_{2^{n}+1}F\right\vert}{2^{n\left( 1/p-2\right)}}
	\end{equation*}
	is bounded from the Hardy space ${{H}_{p}}(G)$ to the Lebesgue space $\text{weak}-L_p(G).$
	
	b) (Sharpness) Let $\varphi :\mathbb{N}\rightarrow \lbrack 1,$ $\infty )$ is a nondecreasing function, satisfying the condition	
	\begin{equation*}
	\overline{\lim_{n\rightarrow \infty }}\frac{2^{n\left( 1/p-1\right)}}{\varphi \left( n\right)}=\infty. 
	\end{equation*}
	Then, there exists sequence $\{f_n, n\in \mathbb{N}\}$  of $p$-atoms, such that	
	\begin{equation*}
	\sup_{n\in \mathbb{N}}\frac{\left\Vert \frac{S_{2^n+1}f_n}{\varphi \left( 2^n+1\right)} \right\Vert_{\text{weak}-L_p(G)}}{\left\Vert f_n \right\Vert_{H_p(G)}}=\infty .
	\end{equation*}
\end{corollary}
Finally, we note that Theorem \ref{theorem1}  implies the following result of Tephnadze \cite{tep3}:

\begin{corollary}
	a) Let $0<p<1$ and $f\in H_p(G).$ Then the weighted  maximal operator $\overset{\sim }{S }^{*,p},$ defined by \eqref{max0}
	is bounded from the martingale Hardy space $H_p(G)$ to the Lebesgue space $\text{weak}-L_p(G).$
	
	b) Let $\{\varphi_n\}$ be any nondecreasing sequence satisfying the condition
	\begin{equation*}
	\overline{\lim_{n\rightarrow \infty }}\frac{\left( n+1\right)^{1/p-1}}{
		\varphi_n}=\infty.
	\end{equation*}
	Then there exists a martingale $f\in H_{p}(G),$ such that
	\begin{equation*}
	\sup_{n\in \mathbb{N}}\left\| \frac{S _{n}f}{\varphi_n }\right\|
	_{p}=\infty .
	\end{equation*}
\end{corollary}

\section{Proofs of the Theorems}

\begin{proof}[\textbf{Proof of Theorem 1.}]
	Since $\sigma _{n}$ is bounded from $L_{\infty }$ to $%
	L_{\infty },$
	by Lemma \ref{W1}, the proof of Theorem \ref{theorem1}
	will be complete, if we prove that	
	\begin{equation}\label{weaktypesigma5}
	t^p\mu\left\{x\in\overline{I_M}:\widetilde{S }^{\ast ,\nabla }a(x) \geq t\right\} \leq c_p<\infty,
	\text{ \ \ \ \ }t\geq 0
	\end{equation}
	for every $p$-atom $a.$ In this parer $c_p$ (or $C_p$) denotes a positive constant depending only on $p$ but which can be different in different places.

	We may assume that $a$ is an arbitrary $p$-atom,
	with support $I, \ \mu \left( I\right) =2^{-M}$ and $I=I_{M}.$ 
	It is easy to see that 
	$S _{n}a\left( x\right) =0, \  \text{ when } \  n< 2^{M}.$
	Therefore, we can suppose that $n\geq 2^{M}.$
	Since $\left\Vert a\right\Vert _{\infty }\leq 2^{M/p},$
	we obtain that 
	\begin{eqnarray*}
		\frac{\left\vert S _{n}a\left( x\right) \right\vert}{2^{\rho\left( n\right) \left( 1/p-1\right) }} 
		&\leq& \frac{1}{2^{\rho\left( n\right) \left( 1/p-1\right)} } 
		\left\Vert a\right\Vert _{\infty }\int_{I_{M}}\left\vert D_{n}\left(
		x+t\right) \right\vert \mu \left( t\right)\\
		&\leq& \frac{1}{2^{\rho\left( n\right) \left( 1/p-1\right)}} 2^{M/p}\int_{I_{M}}\left\vert
		D_{n}\left( x+t\right) \right\vert \mu \left( t\right) .
	\end{eqnarray*}
	
	Let $x\in I_{s}\backslash I_{s+1} ,\,0\leq s< \left[ n\right]
	\leq M$ or $0\leq s\leq M<\left[ n\right].$ Then, it is easy to see that $x+t\in I_{s}\backslash I_{s+1}$ for $t\in I_M$ and if we combine  \eqref{1dn} and \eqref{2dn} we get that 
	$D_{n}\left( x+t\right)
	=0, \  \text{ for } \  t\in I_{M}$ 
	so that 
	\begin{equation}\label{12a}
	\frac{\left\vert S_{n}a\left( x\right) \right\vert}{2^{\rho\left( n\right) \left( 1/p-1\right)}}=0.  
	\end{equation}
	
	Let $I_{s}\backslash I_{s+1} ,\,\left[ n\right] \leq s\leq M$ or $ \left[ n\right] \leq s\leq M.$ Then, it is easy to see that $x+t\in I_{s}\backslash I_{s+1}$ for $t\in I_M$ and if we again combine  \eqref{1dn} and \eqref{2dn} we find that 
	$D_{n}\left( x+t\right)\leq c2^s, \  \text{ for } \  t\in I_{M}$
	and 
	\begin{eqnarray}\label{12}
	\frac{\left\vert S_{n}a\left( x\right) \right\vert }{2^{\rho\left( n\right) \left( 1/p-1\right)}} 
	&\leq& c_p 2^{M/p}\frac{2^{s-M}}{2^{\rho\left( n\right) \left( 1/p-1\right)}} \\ \notag
	&\leq&c_p 
	\frac{2^{[n]\left( 1/p-1\right) +s+M(1/p-1)}}{2^{\left\vert n\right\vert \left( 1/p-1\right)}} \leq c_{p}2^{\left[ n\right] \left( 1/p-1\right) +s}\leq c_{p}2^{s}.
	\end{eqnarray}
	
	By applying (\ref{12a}) and (\ref{12})  for any  $x\in I_{s}\backslash I_{s+1}  ,\,0\leq s< M,$ we find that
	\begin{eqnarray}\label{weaktypesigma1}
		\widetilde{S }^{\ast,\nabla }a\left( x\right)=\sup_{n\in\mathbb{N}}\left(\frac{\left\vert S_{n}a\left( x\right) \right\vert}{2^{\rho\left( n\right) \left( 1/p-1\right)}} \right) \leq C_{p}2^{s/p}.
	\end{eqnarray}
	
	It immediately follows that for  $s\leq M$ we have the following estimate
	\begin{equation*}
	\widetilde{S}^{\ast ,\nabla }a\left( x\right) \leq C_p2^{M/p}\text{ \ \ for any\ \ }x\in I_{s}\backslash I_{s+1}, \ \ s=0,1,\cdots, M
	\end{equation*}%
	and also that
	\begin{equation}\label{weaktypesigma0}
	\mu \left\{ x\in I_{s}\backslash I_{s+1}:\widetilde{S }^{\ast ,\nabla }a\left( x\right)> C_p2^{k/p}\right\} =0, \ \ \ k=M, M+1,\ldots
	\end{equation}

	By combining \eqref{2} and \eqref{weaktypesigma1} we get that
	\begin{eqnarray*}
		\left\{ x\in \overline{I_{N}}:\widetilde{S }^{\ast ,\nabla }a\left( x\right)\geq C_p2^{k/p}\right\} 
		\subset \bigcup_{s=k}^{M-1}\left\{ x\in I_{s}\backslash
		I_{s+1}:\widetilde{S }^{\ast ,\nabla }a\left( x\right)\geq C_p 2^{k/p}\right\} 
	\end{eqnarray*}
	and
\begin{align}\label{weaktypesigma4}
\mu \left\{ x\in \overline{I_{M}}:\widetilde{S }^{\ast ,\nabla }a\left( x\right)\geq C_p 2^{k/p}\right\}\leq  \overset{M-1}{\underset{s=k}{\sum }}\frac{1}{2^{s}}\leq \frac{2}{2^{k}}.
\end{align}
	In view of \eqref{weaktypesigma0} and \eqref{weaktypesigma4} we can conclude that
	\begin{eqnarray*}
		2^{k}\mu \left\{ x\in \overline{I_{N}}:\widetilde{S }^{\ast ,\nabla }a\left( x\right)\geq {C_p}{2^{k/p}}\right\}<c_p<\infty,
	\end{eqnarray*}
which shows that \eqref{weaktypesigma5} holds and the proof of is complete.
\end{proof}

\begin{proof}[\textbf{Proof of Theorem 2.}]
a)	Set
	\begin{equation*}
	f_{n_{k}}\left( x\right) =D_{2^{n_{k}+1}}\left( x\right)
	-D_{2^{{n_{k}}}}\left( x\right) ,\text{ \qquad }n_{k}\geq 3.
	\end{equation*}
	
	It is evident that
	\begin{equation*}
	\widehat{f}_{n_{k}}\left( i\right) =\left\{
	\begin{array}{l}
	1,\ \ \ \text{ if }i=2^{n_{k}},...,2^{n_{k}+1}-1, \\
	0, \ \ \ \text{otherwise}.%
	\end{array}%
	\right.
	\end{equation*}%
	Then we have that
	\begin{equation}\label{14} 
	S_{i}f_{n_{k}}\left( x\right) =\left\{
	\begin{array}{l}
	D_{i}\left( x\right) -D_{2^{n_{k}}}\left( x\right) ,\ \ \text{ if } \ \ 
	i=2^{n_{k}},...,2^{n_{k}+1}-1, \\
	f_{n_{k}}\left( x\right),\ \ \text{ if }\ \ i\geq 2^{n_{k}+1}, \\
	0,\ \ \ \text{ \qquad otherwise}.
	\end{array}%
	\right.  
	\end{equation}
	Since
	\begin{equation}\label{dj}
	D_{j+2^{n_{k}}}\left( x\right) -D_{2^{n_{k}}}\left( x\right)
	=w_{2^{n_{k}}}D_{j}(x),\text{ \qquad }\,j=1,2,..,2^{n_{k}},
	\end{equation}
	from (\ref{1dn}) it follows that
	\begin{eqnarray}\label{15} 
	\left\Vert f_{n_{k}}\right\Vert _{H_{p}}  
	&=&\left\Vert \sup\limits_{n\in \mathbb{N}}S_{2^{n}}f_{n_{k}}
	\right\Vert_p =\left\Vert D_{2^{n_{k}+1}}-D_{2^{n_{k}}}\right\Vert_p\\ \notag
	& =&\left\Vert D_{2^{n_{k}}}\right\Vert_p  \leq 2^{n_{k}(1-1/p)}.  
	\end{eqnarray}
	
	Let  $q_{n_k}^s\in \mathbb{N}$ be such that 
	$2^{n_k}< q_{n_k}^s< 2^{n_k+1},$  where $[q_{n_k}^s]=s,$  $s=0,\cdots, n_k-1.$ 
	By combining (\ref{14}) and \eqref{dj} we can conclude that
	\begin{eqnarray*}
		\left\vert S _{q_{n_k}^s}f_{n_{k}}\left( x\right) \right\vert
		=\left\vert  D_{q_{n_k}^s}\left(
		x\right) -D_{2^{n_{k}}}\left( x\right)  \right\vert =\left\vert  D_{q_{n_k}^s-2^{n_k}-1}\left(x\right) \right\vert
	\end{eqnarray*}

Let $x\in I_{s+1}\left(e_{s}\right)$. By using Lemma \ref{lemma1} we find that
	\begin{equation}\label{jjj}
	\left\vert S_{q_{n_k}^s}f_{n_k}\left( x\right) \right\vert
	\geq c2^{s}
	\end{equation}
	so that
\begin{equation*}
\frac{\left\vert S_{q_{n_k}^s}f_{n_{k}}\left( x\right) \right\vert }{2^{{(1/p-1)}\rho\left( q_{n_k}^s\right) }}\geq \frac{c2^{s/p}}{2^{{n_k}{(1/p-1)} }}.
\end{equation*}
Hence,

\begin{eqnarray}\label{22abc}
	&&\int_{G}\left(\sup_{n\in \mathbb{N}}\frac{\left\vert S_{n}f_{n_{k}}\left( x\right) \right\vert }{2^{{(1/p-1)}\rho\left(n\right)}}\right)^{p}d\mu \left( x\right) \\ \notag
	&\geq &\overset{n_{k}-1}{\underset{s=0}{\sum }}
	\int_{I_{s+1}\left(e_{s}\right)}\left( \frac{\left\vert S _{q_{n^s_{k}}}f_{n_{k}}\left( x\right) \right\vert }{2^{{(1/p-1)}\rho\left(q_{n_k}^{s}\right)} }\right) ^{p}d\mu
	\left( x\right) \\ \notag
	&\geq &c_p\overset{n_{k}-1}{\underset{s=0}{\sum }}\frac{1}{2^{s}}\frac{2^{s}}{2^{n_{k}(1-p)}}\geq \frac{C_p }{2^{n_{k}(1-p)}}\overset{n_{k}}{\underset{s=1}{\sum }}1\geq \frac{C_p n_k }{2^{n_{k}(1-p)}}.
\end{eqnarray}
Finally, by combining (\ref{15}) and \eqref{22abc} we obtain that
\begin{eqnarray*}
	\frac{\left( \int_{G}\left(\sup_{n\in \mathbb{N}}\frac{\left\vert S_{n}f_{n_{k}}\left( x\right) \right\vert }{2^{{(1/p-1)}\rho\left(n\right)}}\right)^{p}d\mu \left( x\right) \right) ^{1/p}}{\left\| f_{n_{k}}\right\| _{H_p}} 
	&\geq&\frac{{\left(\frac{C_pn_k}{2^{n_{k}(1-p)}}\right)}^{1/p}}{2^{n_k(1-1/p)}}\\
	&\geq& c_p n^{1/p}_k\rightarrow \infty,\quad \text{as \quad }k\rightarrow \infty,
\end{eqnarray*}
so the proof of part a) is complete.

b) Under condition \eqref{6aa}  we can choose  $q^{s_k}_{n_k}\in \mathbb{N}$ for some $0\leq s_k <n_k$  such that 
$2^{n_k}< q_{n_k}^{s_k}< 2^{n_k+1}, \ \text{where} \  [q_{n_k}]=s_k$
and
\begin{equation*}
{\lim_{k\rightarrow \infty }}\frac{2^{\rho(q_{n_k}^{s_k}){(1/p-1)}}}{\varphi \left( q_{n_k}^{s_k}\right) }=\infty. 
\end{equation*}

Let $x\in I_{s_k+1}\left(e_{s_k}\right)$. By using  \eqref{jjj} we get that
\begin{equation*}
\left\vert S_{q^{s_k}_{n_k}}f_{n_k}\left( x\right) \right\vert
\geq c2^{s_k} \ \ \ \text{
and } \ \ \ 
\frac{\left\vert S_{q^{s_k}_{n_k}}f_{n_{k}}\left( x\right) \right\vert }{
	\varphi \left( q^{s_k}_{n_k}\right) }\geq \frac{c2^{s_k}}{\varphi \left(q^{s_k}_{n_k}\right) }.
\end{equation*}

Hence, we find that
\begin{eqnarray}\label{16}
&&\mu \left\{ x\in G:\frac{\left\vert S_{q^{s_k}_{n_k}}f_{{n_{k}}}\left( x\right) \right\vert }{\varphi\left(q^{s_k}_{n_k}\right) }\geq \frac{c2^{s_k}}{\varphi
	\left(q^{s_k}_{n_k}\right) }\right\}  \geq \mu \left( I_{s_k+1}(e_{s_k})\right) >c/2^{s_k}.
\end{eqnarray}

By combining (\ref{15}) and (\ref{16}) we get that
\begin{eqnarray*}
	&&\frac{\frac{c2^{s_k}}{\varphi \left( q_{n_k}^{s_k}\right) }\left( \mu \left\{ x\in G:\frac{\left\vert S_{q^{s_k}_{n_k}}f_{{n_{k}}}\left( x\right) \right\vert }{\varphi\left(q_{n_k}^{s_k}\right) }\geq \frac{c2^{s_k}}{\varphi\left(q_{n_k}^{s_k}\right) }\right\} \right) ^{1/p}}{\left\Vert f_{n_{k}} \right\Vert _{H_{p}(G)}} \\
	&\geq &\frac{c_p2^{s_k}}{\varphi \left(q_{n_k}^{s_k}\right)2^{n_k{(1-1/p)}}}\frac{1}{2^{s_k/p}}\\
	&=&\frac{c_p2^{{n_k}{(1/p-1)}}}{2^{{s_k}{(1/p-1)}}\varphi \left(q_{n_k}^{s_k}\right) }=\frac{c_p2^{\rho\left( q_{n_k}^{s_k}\right)(1/p-1)}}{\varphi \left(q_{n_k}^{s_k}\right) }
	\rightarrow \infty \text{,\quad as \quad }k\rightarrow \infty,
\end{eqnarray*}	
so also part b) is proved and the proof is complete.	
\end{proof}

\textbf{Author Contributions: }

DB, LEP and GT  gave the idea and initiated the writing of this paper. SH followed up on this with some complementary ideas.
All authors read and approved the final manuscript. \\
\textbf{Conflict of Interest:}

The authors declare that they have no competing interests. 






\end{document}